\documentclass[10pt]{amsart}
\usepackage{color}

\definecolor{c20}{rgb}{0.,0.7,0.}
\definecolor{c30}{rgb}{0.,0.,1.}
\definecolor{c40}{rgb}{1,0.1,0.7}
\definecolor{c50}{rgb}{1,0,0}
\definecolor{c60}{rgb}{1,0.9,0.1}

\newcommand{\abs}[1]{\left\lvert #1 \right\rvert}

\newcommand{\E}[1]{\mathbb{E}\left(#1\right)}

\newcommand{\pk}[1]{\mathbb{P} \left\{ #1 \right \} }

\newcommand{\R}{\mathbb{R}}

\newcommand{\BQN}{\begin{eqnarray}}
\newcommand{\EQN}{\end{eqnarray}}
\newcommand{\BQNY}{\begin{eqnarray*}}
\newcommand{\EQNY}{\end{eqnarray*}}

\newcommand{\BS}{\begin{sat}}
\newcommand{\ES}{\end{sat}}
\newcommand{\BT}{\begin{theo}}
\newcommand{\ET}{\end{theo}}
\newcommand{\BK}{\begin{korr}}
\newcommand{\EK}{\end{korr}}

\newcommand{\BD}{\begin{de}}
\newcommand{\ED}{\end{de}}

\newcommand{\BIT}{\begin{itemize}}
\newcommand{\EIT}{\end{itemize}}
\newcommand{\BDI}{\begin{description}}
\newcommand{\EDI}{\end{description}}

\newcommand{\BRM}{\begin{remarks}}
\newcommand{\ERM}{\end{remarks}}

\newcommand{\BEL}{\begin{lem}}
\newcommand{\EEL}{\end{lem}}

\newtheorem{theo}{Theorem}[section]
\newtheorem{sat}[theo]{Proposition}
\newtheorem{de}[theo]{Definition}
\newtheorem{lem}[theo]{Lemma}

\newtheorem{example}[theo]{Example}
\newtheorem{korr}[theo]{Corollary}
\newtheorem{remark}[theo]{Remark}
\newtheorem{remarks}[theo]{Remarks}

\newcommand{\nelem}[1]{{Lemma \ref{#1}}}

\newcommand{\netheo}[1]{{Theorem \ref{#1}}}

\newcommand{\prooftheo}[1]{ \textsc{\bf Proof of Theorem} \ref{#1}:}

\newcommand{\prooflem}[1]{\textsc{\bf Proof of Lemma} \ref{#1}:}
\newcommand{\proofkorr}[1]{\textsc{\bf Proof of Corollary} \ref{#1}:}

\newcommand{\COM}[1]{}

\newcommand{\QED}{\hfill $\Box$}

\newcommand{\kb}[1]{\boldsymbol{#1}}
\newcommand{\vk}[1]{\kb{#1}}

\def\vn{\varepsilon}

\def\rw{\rightarrow}

\def\IF{\infty}

\def\invf{\overleftarrow{f}}

\def\LT{\left}
\def\RT{\right}

\def\TT{\mathcal{T}}
\def\oo{(1+o(1))}

\def\HH{\mathcal{H}}

\def\vn{\varepsilon}

\def\to{\rightarrow}

\def\Piter{\mathcal{P}}
\def\MM{\mathcal{M}}
\def\SS{\mathcal S}

\def\DD{\mathcal D}

\def\PTT{ \Upsilon_k (u) }
\def\rd#1{\textcolor{black}{#1}}

\def\rj#1{\textcolor{black}{#1}}

\begin{document}


\title{Tail Asymptotic Behavior of the supremum of a class of chi-square processes}


\author{Lanpeng Ji}
\address{Lanpeng Ji,  School of Mathematics, University of Leeds, Woodhouse Lane, Leeds LS2 9JT, United Kingdom
}
\email{l.ji@leeds.ac.uk}

\author{Peng Liu}
\address{Peng Liu, Department of Statistics and Actuarial Science, University of Waterloo, Canada
}
\email{peng.liu1@uwaterloo.ca}

\author{Stephan Robert}
\address{Stephan Robert, Institute for Information and Communication Technologies, School of Business and Engineering Vaud (HEIG-VD), University of Applied Sciences of Western Switzerland, Switzerland}
\email{Stephan.Robert@heig-vd.ch}

\bigskip

\date{\today}
 \maketitle

 {\bf Abstract:}
We analyze
in this paper the supremum of a class of chi-square processes over non-compact intervals, which can be seen as  a multivariate counterpart of the generalized weighted Kolmogorov-Smirnov statistic.
The boundedness and the exact  tail asymptotic behavior of  the supremum are derived.
 As examples, the
chi-square process generated from the Brownian bridge and the fractional Brownian motion are discussed.



\section{Introduction}
Let $X(t), t\ge 0,$ be a Gaussian process with almost surely (a.s.) continuous sample paths. For a sequence of constants $\{b_i\}_{i=1}^n$ satisfying
$$1=b_1= \cdots= b_k> b_{k+1} \ge \cdots \ge b_n>0$$
we define the chi-square process  as
\BQN\label{eq:chi-p}
\chi_{\vk{b}}^2(t)= \sum_{i=1}^n b_i^2 X_i^2(t), \quad t\ge 0,
 \EQN
  where  $X_i$'s are independent copies of  $X$.
The  supremum of chi-square process   appears naturally as limiting test statistic in various statistical models; see, e.g., \cite{Jar99, JA2, AlbinJA,JAPIT}.  It also plays an important role in reliability applications  in the engineering sciences, see \cite{Lind80, Lindgren1989,HAJI2014, chiLiu} and the references
therein.

Of interest in applied probability and statistics is the tail asymptotics of
$$
 \pk{\sup_{t\in\TT}\chi_{\vk{b}}^2(t)>u},\ \ \ u\to\IF
$$
for an interval $\TT\subset\R_+$, provided that
\BQN \label{eq:chiIF}
 \sup_{t\in\TT}\chi_{\vk{b}}^2(t) <\IF\ \ \ a.s..
\EQN
Numerous contributions have been  devoted to the study of the tail asymptotics of the supremum of chi-square processes over compact intervals $\TT$; see, e.g., \cite{Lindgren1989,konstantinides2004gnedenko,tanH2012, HJi16} and the references therein, where the technique used  is to transform the supremum of  chi-square process into the supremum of a special Gaussian random field.  We refer to, e.g., \cite{AdlerTaylor, Pit96, AZI, MarVad16, CX16, CX16b, BEL17}   for more discussions on the tail asymptotics (or excursion probability) of Gaussian and related fields.

In this paper, we are interested in the analysis of
a class of   weighted locally stationary  chi-square processes  defined by
\BQNY
\sup_{t\in \TT}\frac{\chi^2_{\vk b}(t) }{w^2(t)},  \ \ \mathrm{with \ } \TT=(0,1) \mathrm{\ or\ }  (0,1],  
\EQNY
where $w(\cdot)$ is some positive continuous function definable on the non-compact set $ \TT$, and the generic   process $X$ is the locally stationary Gaussian process.
More precisely,  $X(t), t\in \TT,$ is a centered Gaussian process with a.s. continuous sample paths, unit variance and correlation function $r(\cdot,\cdot)$ satisfying
\BQN\label{cor}
\lim_{h\to 0}\frac{1-r(t,t+h)}{K^2(\abs{h})}=C(t)
\EQN
uniformly in $t\in I$, for all the compact interval $I$ in $\TT$, where $K(\cdot)$ is a positive regularly varying function at 0 with index $\alpha/2\in (0,1]$, and $C(\cdot)$ is a positive continuous function satisfying
$$\lim_{t\to 0}C(t)=\IF \ \ \mathrm{ or\ } \  \lim_{t\to 1}C(t)=\IF.$$
We refer to \cite{LiuJi} for more discussions on  such locally stationary Gaussian processes.

Our motivation for considering the supremum of the  weighted locally stationary  chi-square processes over the non-compact interval $
\TT=(0,1) \ \mathrm{or}\ (0,1]
$ is from its potential applications in statistics.
For instance, in its univariate framework (with $n=1$)
 the following  generalized  weighted Kolmogorov-Smirnov statistic
 \BQNY
W_w := \sup_{t\in(0,1)}  \frac{ \abs{\overline{ B  } (t)}}{ w(t) }, \ \ \text{with} \  \  \overline{B}(t)= \frac{ B(t) }{\sqrt{t(1-t)}},\ \ t\in(0,1),
\EQNY
has been discussed in the statistics literature, see, e.g., 
 \cite{CCHM86},  where $B$ is the standard Brownian bridge with variance function $Var(B(t))=t(1-t), t\in[0,1]$ and
$w$ is a suitably chosen weight function such  that
\BQN\label{eq:Www}
W_w<\IF\ \ \  a.s..
\EQN
 We refer to \cite{CCHM86,DasGupta,  PitTar91, DumKolWil} for further discussions on
 the generalized  weighted Kolmogorov-Smirnov statistic.

An  interesting theoritical question is to find sufficient and necessary conditions on $w$ so that the a.s. finiteness of \eqref{eq:Www} holds. It is shown in \cite{CCHM86}[Theorem 3.3, Theorem 4.2.3] (see also    \cite{DasGupta}[Theorem 26.3]) that
\BQN\label{eq:Ww}
W_w<\IF\ \ a.s. \ \ \Leftrightarrow \ \ \int_0^1 \frac{1}{t(1-t)} e^{-c w^2(t) } dt<\IF  \ \mathrm{for\ some\ } c>0.
\EQN
\COM{
Note that for the normalized standard Brownian bridge
\BQN \label{eq:corrBb}
\lim_{h\to 0}\frac{1-\E{ \overline{B}(t) \overline{B}(t+h)}}{\abs{h}}=C(t)
\EQN
holds uniformly in $t\in I$, any compact interval in $(0,1)$, where $C(t)= \frac{1}{2t(1-t)}$ satisfying $\lim_{t\to 0}C(t)=\lim_{t\to 1}C(t)=\IF$.
Motivated by this we introduce next the class of Gaussian processes that will be discussed.
} 

One of the main results displayed   in \netheo{Thm01law} shows necessary and sufficient conditions on the weight function $w$ under which it holds that
\BQN\label{eq:mXas}
\sup_{t\in \TT}\frac{\chi_{\vk b}^2(t) }{w^2(t)}<\IF\ \ \ \ \ \mathrm{a.s.}.
\EQN
This extends the result of  \eqref{eq:Ww}.
Furthermore, for certain $w$ satisfying \eqref{eq:mXas}  we derive in \netheo{asym} the exact asymptotics of
\BQN\label{eq:mXtail}
\pk{\sup_{t\in \TT}\frac{\chi_{\vk b}^2(t) }{w^2(t)}>u},\   \ \ u\rightarrow\infty.
\EQN
As an important application of \netheo{asym}, we obtain in Corollary 3.4 the tail asymptotics of the supremum of the chi-square process generated from the Brownian bridge. It is worth mentioning that this tail asymptotic result is  new even for the univariate (i.e., $n=1$) case. As a second example, the chi-square process generated by the fractional Brownian motion is  discussed.

We expect that the derived results will have interesting applications in statistics or beyond.


\vskip  0.5 cm
{\it Organization of the rest of the paper}: In Section 2 we  present a preliminary result which is a tailored version of Theorem A.1 in \cite{LiuJi}. The main results are given in Section 3, followed by examples. All the proofs are displayed in Section 4.

\section{Preliminaries}

This section concerns a result  derived in \cite{LiuJi}, which is  crucial for the derivation of \eqref{eq:mXas}.
Based on the discussions therein,  we shall consider  $\int_{0}^{1/2}(C(s))^{1/\alpha}ds=\IF$ or $\int_{1/2}^{1}(C(s))^{1/\alpha}ds=\IF$. 
 For this purpose, of crucial importance is the following function
\BQNY
f(t)=\int_{1/2}^t(C(s))^{1/\alpha}ds,\ \ \ t\in(0,1).
\EQNY
We denote by $\invf(t), t\in(f(0), f(1))$ the inverse function of $f(t),t\in(0,1)$. Further, for any $d>0$,
 let $s_{j,d}^{(1)}=\invf(jd)$, $j\in \mathbb{N}\cup \{0\}$ if $f(1)=\IF$, and let $s_{j,d}^{(0)}=\invf(-jd)$, $j\in \mathbb{N}\cup\{0\}$ if $f(0)=-\IF$. Denote $\Delta_{j,d}^{(1)}=[s_{j-1,d}^{(1)},s_{j,d}^{(1)}], j\in \mathbb{N}$ and
  $\Delta_{j,d}^{(0)}=[s_{j,d}^{(0)},s_{j-1,d}^{(0)}], j\in \mathbb{N}$, which give a partition of $[1/2,1)$ in the case $f(1)=\IF$ and a partition of $(0,1/2]$ in the case $f(0)=-\IF$, respectively. Moreover, let
$q(u)=\overleftarrow{K}(u^{-1/2})$ be the inverse function of $K(\cdot)$ at point $u^{-1/2}$ (assumed to exist asymptotically).

  The following (scenario-dependent) restrictions on
the positive continuous weight function $w^2$ and the correlation function $r(\cdot,\cdot)$ of $X$ play  a crucial role. Let therefore $S\in\{0,1\}$.  \\
\textbf{Condition A}($S$):  The weight function $w^2$ is monotone in a neighbourhood of $S$ and satisfies $\lim_{t\to S}w^2(t)=\IF$. \\
 \textbf{Condition B}$(S)$: Suppose that there exists some constant  $d_0>0$ such that
\BQNY
  \limsup_{j\rw \IF}\sup_{t\neq s\in\Delta^{(S)}_{j,d_0}}\frac{1-r(t,s)}{K^2(|f(t)-f(s)|)} < \IF,
 \EQNY
and when $\alpha=2$ and $k=1$, assume further
 \BQNY
 K^2(|t|)=O(t^2), \ \ t\rw 0.
 \EQNY

\textbf{Condition C}$(S)$: Suppose that there exists some constant  $d_0>0$ such that
\BQNY
    \liminf_{j\rw \IF}\inf_{t\neq s\in\Delta^{(S)}_{j,d_0}}\frac{1-r(t,s)}{K^2(|f(t)-f(s)|)}> 0.
 \EQNY
 Moreover, there exist  $j_0, l_0\in \mathbb{N}$, $M_0, \beta>0, $ such that for $j\geq j_0,$ $l\geq l_0$,
 \BQN\label{B2}
 \sup_{s\in \Delta^{(S)}_{j+l,d_0},t\in \Delta^{(S)}_{j,d_0} }|r(s,t)|<M_0l^{-\beta}.
 \EQN

For the subsequent discussions we present a tailored version of Theorem A.1 of \cite{LiuJi}, focusing on $|f(S)|=\IF$. \rj{We} define
 \BQNY
I_w(S)=\left|\int_{1/2}^S (C(t))^{1/\alpha}\frac{(w(t))^{ k -2}}{q(w^2(t)) }e^{-\frac{w^2(t)}{2}}dt\right|.
 \EQNY

 \BT\label{THM}  Let $X(t), t\in (0,1),$ be a centered locally stationary Gaussian process with a.s. continuous sample paths, unit variance and correlation function $r(\cdot,\cdot)$ satisfying (\ref{cor}) and $r(s,t)<1$ for $s\neq t \in (0,1)$.  Suppose further that, for $S=0$ or $1$, we have $|f(S)|=\IF$ and
   {\bf A}(S), {\bf B}(S) ,{\bf C}(S) are  satisfied. Then
\BQN\label{eq:chig1}
\pk{\chi_{\vk{b}}^2(t)\leq w^2(t) \ \ \text{ultimately as}\ \ t\rw S }=0,\quad \text{or} \quad 1
\EQN
according to
$$I_w(S)=\IF, \quad \text{or}\quad  <\IF.$$

\ET

 \section{Main Results}
In this section, we first give a criteria for \eqref{eq:mXas} to hold  and then display the exact asymptotics of \eqref{eq:mXtail} for different types of $w$ such that \eqref{eq:mXas} is valid.

 \subsection{Analysis of \eqref{eq:mXas}}

Denote by $E(0)=(0,1/2]$ and $E(1)=[1/2,1)$. Recall that $S\in\{0,1\}$. Under the conditions of \netheo{THM}, we have
 that if $I_w(S)<\IF$, then
\BQNY
\sup_{t\in E(S)}\frac{\chi_{\vk b}^2(t) }{w^2(t)}<\IF\ \ \ \ \ \mathrm{a.s.},
\EQNY
however, when $I_w(S)=\IF$  we only see that
\BQNY
\sup_{t\in E(S)}\frac{\chi_{\vk b}^2(t) }{w^2(t)}\ge 1\ \ \ \ \ \mathrm{a.s.}.
\EQNY
Apparently, the above is not informative for the validity of  \eqref{eq:mXas}. On the other hand, it is easily shown that
\BQNY
\sup_{t\in E(S)}\frac{\chi_{\vk b}^2(t) }{w^2(t)}<\IF\ \ a.s.   \  \  \  \Leftrightarrow  \  \  \   \sup_{t\in E(S)}\frac{\abs{X(t)} }{w(t)}<\IF\ \ a.s.,
\EQNY
which means that, instead of the   condition $I_w(S)=\IF$ in \netheo{THM}, a more accurate condition that is independent of  $n,k$ should be possible to ensure that \eqref{eq:mXas} holds.
Inspired by this fact and given the importance of  \eqref{eq:mXas}, we provide below
a sufficient and necessary condition for \eqref{eq:mXas} \rj{to hold}.

Define, for any constant $c>0$ and any positive continuous function $w$
\BQNY
 \quad J_{c, w}(S)=\left|\int_{1/2}^S (C(t))^{1/\alpha}e^{-c w^2(t)}dt\right|.
 \EQNY

Below is our first principal result, a criterion for \eqref{eq:mXas},  which is a generalization of \eqref{eq:Ww}.

\BT\label{Thm01law}
Under the conditions of \netheo{THM} we have 
\BQNY
\sup_{t\in E(S)}\frac{\chi_{\vk b}^2(t) }{w^2(t)}<\IF\ \ a.s.   \  \  \  \Leftrightarrow  \  \  \  J_{c,w}(S)<\IF \quad \text{for some}\quad c>0.
\EQNY

\ET

\COM{
\begin{remark}
It is not surprising that $J_{c, w}(S)$ appeared above does not depend on $n,k$ (compare with $I_w(S)$); this is due to the fact that 
\BQNY
\sup_{t\in E(S)}\frac{\chi_{\vk b}(t) }{w(t)}<\IF\ \ a.s.   \  \  \  \Leftrightarrow  \  \  \   \sup_{t\in E(S)}\frac{\abs{X(t)} }{w(t)}<\IF\ \ a.s.
\EQNY

\end{remark}
}

Next we illustrate the criteria presented in 
Theorem \ref{Thm01law} by  an example of a weighted chi-square process with generic process being the normalized standard Brownian bridge, which further provides us with a clear comparison between $I_w(S)$ and $J_{c,w}(S)$.

\begin{example}\label{example}
Let $X(t)=\overline{B}(t), t\in(0,1)$, and, with $\rho_1>0, \rho_2\in \mathbb{R}$, define
 \BQN\label{wrho}
 w^2_{\rho_1,\rho_2}(t)=2\rho_1\ln\ln\LT(\frac{e^2}{t(1-t)}\RT)+2\rho_2\ln\ln\ln\LT(\frac{e^2}{t(1-t)}\RT),\ t\in(0,1).
 \EQN
First note that for the normalized standard Brownian bridge
\BQN \label{eq:corrBb}
\lim_{h\to 0}\frac{1-\E{ \overline{B}(t) \overline{B}(t+h)}}{\abs{h}}= \frac{1}{2t(1-t)}
\EQN
holds uniformly in $t\in I$, for any compact interval $I$ in $(0,1)$.
This means that \rj{$ \overline{B}$ is a locally stationary Gaussian process with}
 $$
 K(h)=\sqrt{\abs{h}},\ \ \alpha=1, \ \ q(u)=u^{-1}.
 $$
 Furthermore,
 $$
 f(t)=\int_{1/2}^t \frac{1}{2s (1-s)} ds =\frac{1}{2}\ln\LT(\frac{t}{1-t} \RT)
 $$
 implying that $f(1)=-f(0)=\IF$. Moreover,  by the proof of Corollary 2.6 in \cite{LiuJi} we have that
conditions {\bf B(S)} and {\bf C(S)} are satisfied by $\overline{B}(t), t\in (0,1)$, and $\E{\overline{B}(t), \overline{B}(s)}<1$ for $s\neq t, s,t\in (0,1)$. Thus, all the conditions of \netheo{Thm01law} are fulfilled.

\rj{Next, on one} hand, we have
 \BQNY
 \frac{1}{t(1-t)}(w_{\rho_1,\rho_2}(t))^{k}e^{-\frac{w^2_{\rho_1,\rho_2}(t)}{2}}\sim \frac{Q}{t(1-t)\LT(\ln\LT(\frac{1}{t(1-t)}\RT)\RT)^{\rho_1}\LT(\ln\ln \LT(\frac{e^2}{t(1-t)}\RT)\RT)^{\rho_2-k/2}}
 \EQNY
as $t\rw 0$ or $t\rw 1$,
 with $Q$  some positive constant. Thus, elementary calculations show that
\BQNY
I_w(0)=I_w(1)= \int_{1/2}^1\frac{(w_{\rho_1,\rho_2}(t))^{k }}{t(1-t)}e^{-\frac{w^2_{\rho_1,\rho_2}(t)}{2}}dt<\IF 
\EQNY
holds if and only if
\BQN\label{eq:rho}
\rho_1>1, \quad \text{or\ }  \rho_1=1 \  \text{and} \  \rho_2>1+k/2.
\EQN
On the other hand, we can show  that the functions $w_{\rho_1,\rho_2}(t)$ satisfying that  $\exists c>0$ such that $J_{c, w}(S)<\IF$ are not restricted to the ones satisfying \eqref{eq:rho}.
In fact, since for any $\rho_1>0$ there exists some $c$ such that  $\rho_1> \frac{1}{2c}$,
we have that
\BQNY
J_{c,w}(0)=J_{c,w}(1)&=& \int_{1/2}^1\frac{1}{t(1-t)}e^{- cw^2_{\rho_1,\rho_2}(t) }dt\\
&\le& \int_{1/2}^1\frac{1}{t(1-t) \LT(\ln\LT(\frac{1}{t(1-t)}\RT)\RT)^{2c\rho_1}\LT(\ln\ln \LT(\frac{e^2}{t(1-t)}\RT)\RT)^{2c\rho_2}}dt<\IF
\EQNY
holds  for any 
$\rho_2\in \mathbb{R}$.
Thus, we conclude from  \netheo{Thm01law} that
\BQNY
\sup_{t\in (0,1)}\frac{\chi_{\vk b}^2(t) }{w_{\rho_1, \rho_2}^2(t)}<\IF\ \ a.s.
\EQNY
holds for any $\rho_1>0$ and $\rho_2\in \mathbb{R}$.

\end{example}
The exact tail asymptotics of $\sup_{t\in (0,1)}\frac{\chi_{\vk b}^2(t) }{w_{\rho_1,\rho_2}^2(t)}$ will be discussed in next section.

 \subsection{Asymptotics of \eqref{eq:mXtail}}

For those $w$ such that \eqref{eq:mXas} holds, of interest is the exact tail asymptotic behavior of $\sup_{t\in \TT} \frac{\chi_{\vk b}^2(t) }{w^2(t)}$. Actually, as we have seen,  the behavior of $w$ around 0 and 1 plays a crucial  role for the finiteness in \eqref{eq:mXas}. However, this  does not apply to the tail asymptotics of $\sup_{t\in \TT} \frac{\chi_{\vk b}^2(t) }{w^2(t)}$. It turns out that only the probability mass in the neighborhood of  minimizer of $w$ contribute to the tail asymptotics, indicating that the other part of the process including the part around 0 or 1 can be neglected. As discussed in \cite{PhysRevE.86.041115}, the weight function is introduced when constructing the Goodness-of-Fit test which is intended to emphasize a specific region of the domain. 
With these motivations,  for the tail asymptotics we shall consider the following two types of $w$:

{\bf Assumption F1:} The function $w$ attains its minimum  at finite distinct inner points $\{t_i\}_{ i=1}^ m$ of $\TT$, and
\BQN\label{eq:g}
w(t_i+t)=w(t_i)+a_i \abs{t_i}^{\beta_i}\oo,\  \ \  t\to t_i 
\EQN
holds for some positive constants $a_i,\beta_i>0, i=1,2,\dots, m$.

{\bf Assumption F2:} The function $w$ attains its minimum at all points on disjoint intervals $[c_i,d_i]\subseteq\TT$, $i=1,2,\dots, m$ (i.e., $w$ is a constant on these intervals).

Under assumption {\bf F1}, we need additional conditions which are stated below.
Recall
$q(u)=\overleftarrow{K}(u^{-1/2})$. It follows that $q(u)$ is a regularly \rj{varying} function at infinity with index $-1/\alpha$ which can be further expressed as $q(u)=u^{-1/\alpha}L(u^{-1/2})$, with
$L(\cdot)$  a slowly varying function at $0$.   Denote further $\beta=\max_{1\le i\le m}\beta_i$. According  to the values of $L(u^{-1/2})$ as $u\to\IF$, we consider the following three scenarios:\\ 
\textbf{C1($\beta$):} $\beta>\alpha$, or $\beta=\alpha$ and $\lim_{u\rw\IF} L(u^{-1/2})=0$;\\
\textbf{C2($\beta$):} $\beta=\alpha$ and $\lim_{u\rw\IF }L(u^{-1/2})=\mathcal{L}\in(0,\IF)$;\\
\textbf{C3($\beta$):} $\beta<\alpha$, or $\beta=\alpha$ and $\lim_{u\rw\IF}L(u^{-1/2})=\IF.$

Before displaying our results, we introduce two important constants. One is the   {\it Pickands constant}     defined by
\BQNY\label{pick}
\mathcal{H}_{2H}=\lim_{S\rightarrow\infty}\frac{1}{S}\E{ \exp\biggl(\sup_{t\in[0,S]}\Bigl(\sqrt{2}B_H(t)-t^{2H}\Bigr)\biggr)},
\EQNY
with $B_H(t),t\in\R,$ a standard  fractional Brownian motion (fBm)  defined on $\R$ with  Hurst index  $H\in (0,1]$. And the other one is the
  {\it Piterbarg  constant}  defined by
\BQNY
    \Piter_{2H}^{d}=\lim_{\lambda\to\IF}   \E{\exp\left(\sup_{t\in[-\lambda,\lambda ]}\left(\sqrt{2}B_{H}(t)-(1+d)\abs{t}^{2H}\right)\right)},\ \ \ d>0.
\EQNY

 We refer to \cite{Pit96, debicki2002ruin, DeMan03, DikerY, Bai17, DH17} for the properties and generalizations of the Pickands-Piterbarg type constants. In what follows, $\alpha$ will play a similar role as $2H.$
 Moreover, We shall use the standard   notation for asymptotic equivalence of two functions $f$ and $h$. Specifically, we  write
$f(x)\sim h(x)$, if  $ \lim_{x \to a}  {f(x)}/{h(x)} = 1$  ($a\in\R\cup\{\IF\}$), and further, write $ f(x) = o(h(x))$, if $ \lim_{x \to a}  {f(x)}/{h(x)} = 0$.

Let  $K=\{1\leq i\leq m: \beta_i=\beta \}$ and $K^c=\{ 1\leq i\leq m: \beta_i<\beta\}$. Below is our second principal result.
\BT\label{asym}
 Let $ \frac{\chi_{\vk b}^2(t) }{w^2(t)}, {t\in \TT},$ be the weighted locally stationary chi-square process considered in \netheo{THM} such that \eqref{eq:mXas} holds. We have:

(i). If {\bf F1} is satisfied, then, as $u\rw\IF$,
\BQNY
\pk{\sup_{t\in \TT} \frac{\chi_{\vk b}^2(t) }{w^2(t)}>u}  \sim \left(\prod_{i=k+1}^n (1-b_i^2)^{-1/2}\right) \MM (u ) \ \Upsilon_k( w^2(t_1)u),
\EQNY
where (with the convention $\prod_{i=p}^q=1$ if $q<p$)
\BQN\label{eq:UU}
 \PTT:=\pk{ \chi_{k,\vk{1}}^2(0) >u}= \frac{2^{(2-k)/2}}{\Gamma(k/2)}u^{k/2-1}\exp\left(-\frac{u}{2}\right),\ \ u>0,
\EQN
and
\BQNY\label{eq:MM}
 \MM ( u)= \left\{
          \begin{array}{ll}
  2\left(\sum_{i\in K}a_i^{-1/\beta} (C(t_i))^{1/\alpha}\right) \rd{(w(t_1))^{2/\alpha-1/\beta}} \Gamma(1/\beta+1)  \mathcal{H}_{\alpha}  (q  (u))^{-1}  u^{ - {1}/{\beta} }, & \hbox{for} \  \textbf{C1($\beta$)} ,\\
\sum_{i\in K} \Piter_{\alpha}^{  a_i(w(t_1)C(t_i))^{-1}\mathcal{L}^\alpha }+ \sharp K^c,&  \hbox{for} \  \textbf{C2($\beta$)},\\
m, &  \hbox{for} \  \textbf{C3($\beta$)}.
              \end{array}
            \right.
\EQNY

(ii). If {\bf F2} is satisfied, then, as $u\rw\IF$,
\BQNY
\pk{\sup_{t\in \TT} \frac{\chi_{\vk b}^2(t) }{w^2(t)}>u}  \sim  \left(\prod_{i=k+1}^n\left(1-b_i^2\right)^{-1/2}\right)   \left( \sum_{j=1}^m\int_{c_j}^{d_j} (C(t))^{1/\alpha}dt\right)  \mathcal{H}_{\alpha}\ (q(w^2(c_1)u) )^{-1}  \Upsilon_k( w^2(c_1)u).
\EQNY
\ET

We conclude this section with two  applications of \netheo{asym}.
The first one is on the weighted locally stationary chi-square process discussed in Example \ref{example}, and the second one concerns the weighted locally stationary chi-square process with generic  process $X$ being a normalized standard fBm.

\BK\label{corollary} Let $\frac{\chi_{\vk b}^2(t) }{w_{\rho_1,\rho_2}^2(t)}, t\in(0,1)$, with  $\rho_1>0$ and $\rho_2\in\mathbb{R}$, be the weighted locally stationary chi-square process discussed in Example \ref{example}. We have, as $u\to\IF,$
if $ \rho_2\ge  -\rho_1\ln\ln(4e^2)$, then
\BQNY
\pk{\sup_{t\in (0,1)}\frac{\chi_{\vk b}^2(t) }{w_{\rho_1, \rho_2}^2(t)}>u}\sim \left(\prod_{i=k+1}^n (1-b_i^2)^{-1/2}\right)  \MM(u)\Upsilon_k( 2A_1 u),
\EQNY
where $A_1= \rho_1\ln\ln\LT(4e^2 \RT)+ \rho_2\ln\ln\ln\LT(4e^2 \RT)$ and
\BQNY
 \MM ( u)= \left\{
          \begin{array}{ll}
2A_1\sqrt{\frac{\pi\ln\LT(4e^2 \RT)\ln\ln\LT(4e^2 \RT)}{\rho_1 \ln\ln\LT(4e^2 \RT)+\rho_2}}  u^{1/2}, &  \hbox{for} \  \rho_2>  -\rho_1\ln\ln(4e^2),\\
2\Gamma(1/4) A_1 \LT(\frac{  \ln\ln\LT(4e^2 \RT) \LT(\ln(4e^2 )\RT)^2}{\rd{8\rho_1} }\RT)^{1/4}  u^{3/4},&  \hbox{for} \  \rho_2=  -\rho_1\ln\ln(4e^2),\\
              \end{array}
            \right.
\EQNY
and if $ \rho_2<  -\rho_1\ln\ln(4e^2)$, then
\BQNY
\pk{\sup_{t\in (0,1)}\frac{\chi_{\vk b}^2(t) }{w_{\rho_1, \rho_2}^2(t)}>u}\sim  \rd{2A_2} \left(\prod_{i=k+1}^n (1-b_i^2)^{-1/2}\right) \rho_1^{-1} Q  \sqrt{ -2\pi \rho_2 }  u^{1/2}\Upsilon_k( 2A_2 u),
\EQNY
where  $A_2=\rho_2(\ln(-\rho_2) -\ln(\rho_1)-1)$ and
$$
Q=  \rd{ \frac{1 }{2t_1-1} \ln\LT(\frac{e^2}{t_1(1-t_1)}\RT)},\ \ \
t_1= 1/2+\sqrt{1/4- e^{2-e^{-\rho_2/\rho_1}}} .
$$
\EK

Next, we consider $B_H(t), t\ge 0,$ to be the standard fBm with Hurst index $H\in (0,1)$ and covariance function
$$
Cov(B_H(s),B_H(t))=\frac{1}{2}\LT(|s|^{2H}+|t|^{2H}-|s-t|^{2H}\RT), \quad s,t\geq 0.
$$
Denote by $\overline{B}_H(t)={B_H(t)}/{t^H}, t\in(0,1]$ the normalized standard fBm defined on $(0,1]$.
Further, for any $\rho>0$ and $\vn\in(0,1)$, we define
\BQN\label{eq:rv}
w_{\rho,\vn}^2 ( t)= \left\{
          \begin{array}{ll}
          \rho\ln\ln\LT(e^2/t\RT), &  \hbox{for} \  t\in(0,\vn),\\
\rho\ln\ln\LT(e^2/\vn\RT),&  \hbox{for} \ t\in[\vn,1].
              \end{array}
            \right.
\EQN
We have the following result.

\BK\label{corollary1}
 Let $\frac{\chi_{\vk b}^2(t) }{w_{\rho,\vn}^2(t)}, t\in(0,1],$ be a weighted chi-square process with generic process $\overline{B}_H(t), t\in (0,1]$ and $w_{\rho,\vn}^2$ given in \eqref{eq:rv}.
 Then, we have, as $u\to\IF,$
\BQNY
\pk{\sup_{t\in (0,1)}\frac{\chi_{\vk b}^2(t) }{w_{\rho, \vn}^2(t)}>u} &\sim & \left(\prod_{i=k+1}^n (1-b_i^2)^{-1/2}\right)   (-\ln(\vn)) \LT(\ln\ln\LT(e^2/\vn\RT) \rho/2 \RT)^{\frac{1}{2H}}   \\
&& \times \mathcal{H}_{2H} u^{\frac{1}{2H}}\Upsilon_k( \rho \ln\ln\LT(e^2/\vn\RT)  u).
\EQNY
\EK

\section{Proofs}

This section is devoted to the proof of all the results presented in Section 3.

\prooftheo{Thm01law}
\COM{Let $\widetilde{I}_g(S)$ denote
$$\widetilde{I}_g(S)=\left|\int_{1/2}^S (C(t))^{1/\alpha}\frac{(g(t))^{-1}}{q(g^2(t)) }e^{-\frac{g^2(t)}{2}}dt\right|.$$
 Then it follows from Theorem 5.1 in \cite{LiuJi} that
 \BQN\label{eq1}
 \widetilde{I}_{Mg}(S)<\IF \quad \Rightarrow \quad \pk{X^2(t)\leq M^2g^2(t) \ \ \text{ultimately as}\ \ t\rw S }=1
 \EQN
 holds if {\bf B(S)} is satisfied,
 and if {\bf C(S)} is satisfied,
\BQN\label{eq2}
\widetilde{I}_{Mg}(S)=\IF \quad \Rightarrow \quad \pk{X^2(t)\leq M^2g^2(t) \ \ \text{ultimately as}\ \ t\rw S }=0.
\EQN
}
Note that $t^{k/2-1}(q(t))^{-1}$ is a positive regularly varying function at $\IF$ with index $\kappa=k/2-1+1/\alpha\ge 0$. Thus, by Potter bound (e.g., \cite{BI1989})
$$c_1t^{\kappa-1}\leq t^{k/2-1}(q(t))^{-1}\leq c_2t^{\kappa+1}, \ \ t\geq c_3,$$
holds for some constants $c_1, c_2, c_3>0$ ,
which, together with the fact that $w^2(t)\rw \IF$ as $t\rw S$, leads to
\BQN\label{eq3}
 Q_1e^{-w^2(t)}\leq \frac{(w(t))^{ k -2}}{q(w^2(t)) } e^{-\frac{w^2(t)}{2}}\leq  Q_2e^{-\frac{w^2(t)}{3}}
\EQN
for all $t$ approaching $S,$ with some positive constants $Q_1,Q_2.$
Therefore, if
$J_{c, w}(S)<\IF$ holds for some $c>0$, then, by \eqref{eq3},
$$I_{\sqrt{3c} w}(S)<\IF.$$
This together with  ii) of \netheo{THM} yields that
$$\limsup_{t\rw S}   \frac{\chi_{\vk b}^2(t) }{w^2(t)} \leq  3c  \ \ a.s.$$
showing that
$$\sup_{t\in E(S)}  \frac{\chi_{\vk b}^2(t) }{w^2(t)} <\IF\ \ a.s.$$
On the other hand, if $J_{c,w}(S)=\IF$ for all $c>0$, then, by \eqref{eq3},
$I_{\sqrt c w}(S)=\IF$. Thus, by  iii) of \netheo{THM}
$$ \limsup_{t\rw S }  \frac{\chi_{\vk b}^2(t) }{w^2(t)}  \geq   c  \ \ a.s. $$
holds for all $c>0,$ implying that
$$ \sup_{t\in E(S)}  \frac{\chi_{\vk b}^2(t) }{w^2(t)}  =\IF\ \ a.s.$$
This completes the proof. \QED

We show next a  version of the  Borell-TIS inequality for chi-square process, which will play a key role in the proof of \netheo{asym}. We refer to, e.g., \cite{AdlerTaylor, GennaBorell}  for discussions on the Borell-TIS inequality for Gaussian random fields. Denote below $\SS\subseteq \R$ to be any fixed interval.
\BEL\label{ThmBorell}
Let $\chi_{\vk{b}}^2(t), t\in \SS,$ be a chi-square process with generic  centered Gaussian process $X$ which has a.s. continuous sample paths and variance function denoted by $\sigma_X^2(t)$. If
$$
\sup_{t\in\SS} X(t)<\IF \ \ \ \ a.s.,
$$
then there exists some positive constant $Q$ such that for all $u>Q^2$ we have
\BQN\label{eq:Borell}
\pk{\sup_{t\in \SS} \chi_{\vk b}^2(t) >u}\le \exp\LT(-\frac{(\sqrt u-Q)^2}{2\sup_{t\in \SS}\sigma^2_X(t)}\RT).
\EQN

\EEL

\prooflem{ThmBorell} Using the classical approach when dealing with  chi-square processes as, e.g., in \cite{JAPIT, Pit96, LiuJi}, we introduce a particular Gaussian random field, namely,
\BQN\label{Y}
Y_{\vk{b}}(t,\vk{\theta}):=\sum_{i=1}^nb_iX_i(t)v_i(\vk{\theta}), \ \ (t,\vk{\theta})\in \DD=:\SS\times[-\pi,\pi]\times[-\pi/2,\pi/2]^{n-2},
\EQN
where $\vk{\theta}=(\theta_2,\theta_3,\cdots,\theta_n),$ and $v_n(\vk{\theta})=\sin( \theta_n), v_{n-1}(\vk{\theta})=\sin (\theta_{n-1})\cos(\theta_n),\cdots, v_1(\vk{\theta})=\cos (\theta_n)\cdots\cos( \theta_2)$ are spherical coordinates. In view of  \cite{Pit96},  for any $u>0$
\BQN\label{eq:Sph}
\mathbb{P}\left(\sup_{t\in \SS}  \chi_{\vk{b}}^2 (t)  >u\right)=\mathbb{P}\left(\sup_{(t,\vk{\theta})\in  \DD }Y_{\vk{b}} (t,\vk{\theta})> \sqrt u \right).
\EQN
 Since  the variance function of $Y_{\vk{b}}$ satisfies for $u>0$
\BQNY
\E{\left(Y_{\vk{b}}(t,\vk{\theta})\right)^2}  =  \sigma_X^2(t)\left(1-(1-b_n^2)\sin^2(\theta_n)-\sum_{i=k+1}^{n-1}
(1-b_i^2)\left(\prod_{j=i+1}^{n}\cos^2(\theta_j)\right)\sin^2 (\theta_i)\right),
\EQNY
we have
\BQN\label{upperbound}
\sup_{(t,\vk{\theta})\in  \DD } \E{\left(Y_{\vk{b}}(t,\vk{\theta})\right)^2}  \le  \sup_{t\in\SS} \sigma_X^2(t).
\EQN
Then, by \eqref{eq:Sph} and the Borell-TIS inequality for Gaussian random fields (cf.  \cite{AdlerTaylor}[Theorem 2.1.1]) we conclude that \eqref{eq:Borell} holds with $Q=\E{\sup_{(t,\vk{\theta})\in  \DD }Y_{\vk{b}} (t,\vk{\theta})}<\IF$. This completes the proof. \QED

The next result concerns a upper bound for the tails of double-sup of the locally stationary chi-square processes, which will also play a key role in the proof of \netheo{asym}.

\BEL\label{LemDouble}
Let $\chi_{\vk{b}}^2(t), t\in \SS,$ be a  chi-square process with the generic  centered locally stationary Gaussian process $X$ which has a.s. continuous sample paths. If further the correlation function of $X$ satisfies
\BQN\label{eq:rr}
r(s,t)<1  \ \mathrm{for\ any\ } s\neq t\in\SS,
\EQN
then,  for any compact intervals $\SS_1,\SS_2\subset \SS$ such that $\SS_1\cap \SS_2=\emptyset$ we have
\BQNY 
\pk{\sup_{t\in \SS_1} \chi_{\vk b}^2(t) >u, \ \sup_{t\in \SS_2} \chi_{\vk b}^2(t) >u}\le \exp\LT(-\frac{(2\sqrt u-Q)^2}{2(2+2\eta)}\RT).
\EQNY
for all $u>Q^2$, with some constant $Q>0$ and $\eta\in(0,1)$.
\EEL

\prooflem{LemDouble} Using the expression of $Y_{\vk{b}}(t,\vk{\theta})$ given in (\ref{Y}), we have
\BQNY
\pk{\sup_{t\in \SS_1} \chi_{\vk b}^2(t) >u, \ \sup_{t\in \SS_2} \chi_{\vk b}^2(t) >u}
&= & \pk{\sup_{(t,\vk{\theta})\in  \DD_1}Y_{\vk{b}} (t,\vk{\theta})> \sqrt{u}, \sup_{(t,\vk{\theta})\in  \DD_2 }Y_{\vk{b}} (t,\vk{\theta})> \sqrt{u}}\\
&\leq&  \pk{\sup_{(t,\vk{\theta})\in  \DD_1, (t',\vk{\theta}')\in  \DD_2 } (Y_{\vk{b}} (t,\vk{\theta})+Y_{\vk{b}} (t',\vk{\theta}'))>2 \sqrt{u}}
\EQNY
where
$$\DD_i=\SS_i\times[-\pi,\pi]\times[-\pi/2,\pi/2]^{n-2},\quad  i=1,2.$$
By \eqref{eq:rr} we have that there exists some $\eta\in(0,1)$ such that
\BQN\label{doubleupper}
\E{\left(Y_{\vk{b}} (t,\vk{\theta})+Y_{\vk{b}} (t',\vk{\theta}')\right)^2}&=&\E{\left(Y_{\vk{b}}(t,\vk{\theta})\right)^2}+\E{\left(Y_{\vk{b}}(t',\vk{\theta}')\right)^2}+2\sum_{i=1}^n \E{X_i(t)X_i(t')}b_i^2v_i(\vk{\theta})v_i(\vk{\theta}')\nonumber\\
&\leq& 2+2\eta\sum_{i=1}^n b_i^2v_i(\vk{\theta})v_i(\vk{\theta}')\nonumber\\
&\leq& 2+2\eta, \quad (t,\vk{\theta})\in  \DD_1, (t',\vk{\theta}')\in  \DD_2.
\EQN
Consequently, by the Borell-TIS inequality
\BQNY
 \pk{\sup_{(t,\vk{\theta})\in  \DD_1, (t',\vk{\theta}')\in  \DD_2 } (Y_{\vk{b}} (t,\vk{\theta})+Y_{\vk{b}} (t',\vk{\theta}'))>2 \sqrt{u}} \leq
\exp\LT(-\frac{(2  \sqrt u-Q)^2}{2(2+2\eta)}\RT)
\EQNY
for all \rd{$u>\frac{Q^2}{4}$, with $Q=\E{\sup_{(t,\vk{\theta})\in  D_1, (t',\vk{\theta}')\in  D_2 } (Y_{\vk{b}} (t,\vk{\theta})+Y_{\vk{b}} (t',\vk{\theta}'))}<\IF$}. Thus, the claim follows. This completes the proof.
\QED

\prooftheo{asym} Without loss of generality, we show the proof only for the case where $\TT=(0,1)$.

\underline{(i).} Let $\rho>0$ be a sufficiently small constant such that
$$
[t_i-\rho,t_i+\rho]\cap [t_j-\rho,t_j+\rho]=\emptyset,\ \ \ \mathrm{for\ all\ } i\neq j.
$$
 Further, denote $\TT_\rho=\TT\setminus \bigcup_{i=1}^m[t_i-\rho,t_i+\rho]$.
 It follows from the Bonferroni inequality (e.g., \cite{Michna09})  that
\BQN\label{decom}
&&\sum_{i=1}^m p_i(u)+\pk{\sup_{t\in \TT_\rho} \frac{\chi_{\vk b}^2(t) }{w^2(t)}>u}\nonumber \\
&&\geq  \pk{\sup_{t\in \TT} \frac{\chi_{\vk b}^2(t) }{w^2(t)}>u}\\
&&\geq \sum_{i=1}^m p_i(u)-\sum_{1\leq i<j\leq m}\pk{\sup_{t\in [t_i-\rho,t_i+\rho]} \frac{\chi_{\vk b}^2(t) }{w^2(t)}>u, \sup_{t\in [t_j-\rho,t_j+\rho]} \frac{\chi_{\vk b}^2(t) }{w^2(t)}>u}, \nonumber
\EQN
where
$$p_i(u)=\pk{\sup_{t\in [t_i-\rho,t_i+\rho]} \frac{\chi_{\vk b}^2(t) }{w^2(t)}>u}.$$
We first focus on the asymptotics of $p_i(u)$ as $u\to\IF$. Denote
\BQNY
Y(t)=\frac{w(t_1)}{w(t)}X(t), \ \ t\in\TT.
\EQNY
We have
$$
p_i(u)=\pk{\sup_{t\in [t_i-\rho,t_i+\rho]} \sum_{l=1}^nb_l^2 Y_l^2(t)  >w^2(t_1) u}, \quad 1\leq i\leq m,
$$
where $\{Y_l\}_{i=1}^n$ is a sequence of independent copies of Gaussian process $Y.$
It can be shown that, by {\bf F1}, for any $i=1,2,\cdots, m$,
$$
\sigma_Y(t)=\sqrt{\E{(Y(t))^2}}=\frac{w(t_1)}{w(t)}, \ \ t\in [t_i-\rho,t_i+\rho],
$$
attains its  maximum which is equal to 1 at the unique  point  $t_i$, and further
\BQNY
\sigma_Y(t_i+t)=1-\frac{a_i}{w(t_1)}\abs{t}^{\beta_i}\oo,\ \ t\to 0.
\EQNY
Moreover, by \eqref{cor}
\BQNY
1-Corr\left(Y(t_i+t), Y(t_i+s)\right)=C(t_i)K^2(|t-s|)\oo,\ \ t\to 0.
\EQNY
Consequently, it follows from \cite{chiLiu}[Theorem 5.2]  that, as $u\to\IF$,
\BQN\label{eq:pu}
p_i(u) \sim\left( \prod_{l=k+1}^n (1-b_l^2)^{-1/2}\right) \MM _i(\beta_i, u ) \ \Upsilon_k( w^2(t_1)u),
\EQN
 where $\Upsilon_k( \cdot)$ is given in \eqref{eq:UU} and
\BQNY
 \MM_i (\beta_i, u)= \left\{
          \begin{array}{ll}
  2a_i^{-1/\beta_i}\rd{(w(t_1))^{2/\alpha-1/\beta_i}}(C(t_i))^{1/\alpha} \Gamma(1/\beta_i+1)  \mathcal{H}_{\alpha}  (q  (u))^{-1}  u^{ - {1}/{\beta_i} }, & \hbox{for} \  \textbf{C1($\beta_i$)} ,\\
 \Piter_{\alpha}^{  a_i(w(t_1)C(t_i))^{-1}\mathcal{L}^\alpha } &  \hbox{for} \  \textbf{C2($\beta_i$)},\\
1, &  \hbox{for} \  \textbf{C3($\beta_i$)}.
              \end{array}
            \right.
\EQNY
In the sequel, we discuss the three scenarios {\bf C1($\beta$), C2($\beta$), C3($\beta$)}   one-by one.

{\it \underline{\bf C1($\beta$)}}.    Using the fact that $\beta=\max_{i=1}^m \beta_i$, we have that
\BQNY
  \MM_j (\beta_j, u)=o\left(\MM_i (\beta_i, u)\right), \quad u\rw\IF
\EQNY
for any $i\in K$ and $j\in K^c$. This implies that
\BQN\label{main}
\sum_{i=1}^m p_i(u)\sim \sum_{i\in K} p_i(u)\sim \left( \prod_{l=k+1}^n (1-b_l^2)^{-1/2}\right) \MM ( u ) \ \Upsilon_k( w^2(t_1)u),
\EQN
where $\MM(\cdot)$ is given in \eqref{eq:MM}.
On the other hand, it follows directly from \nelem{ThmBorell} that
\BQNY
\pk{\sup_{t\in \TT_\rho} \frac{\chi_{\vk b}^2(t) }{w^2(t)}>u}\le  \exp\LT(-\frac{\inf_{t\in \TT_\rho}w^2(t) (\sqrt u-Q)^2}{2}\RT)
\EQNY
holds for all $u>Q^2$, with $Q$ some positive constant. Since further, by {\bf F1},
$$
\inf_{t\in \TT_\rho}w^2(t)>w^2(t_1),
$$
we have that
\BQN\label{eq:TTr}
\pk{\sup_{t\in \TT_\rho} \frac{\chi_{\vk b}^2(t) }{w^2(t)}>u}=o\left(\MM ( u ) \ \Upsilon_k( w^2(t_1)u)\right),\ \ u\to\IF.
\EQN
Moreover, since for any $i\neq j$
\BQNY
&& \pk{\sup_{t\in [t_i-\rho,t_i+\rho]} \frac{\chi_{\vk b}^2(t) }{w^2(t)}>u, \sup_{t\in [t_j-\rho,t_j+\rho]} \frac{\chi_{\vk b}^2(t) }{w^2(t)}>u}\\
&&\quad \le  \pk{\sup_{t\in [t_i-\rho,t_i+\rho]}  \chi_{\vk b}^2(t)  >w^2(t_1) u, \sup_{t\in [t_j-\rho,t_j+\rho]}  \chi_{\vk b}^2(t)  > w^2(t_1)u}.
\EQNY

we have from \nelem{LemDouble} that, for all $u$ large,
\BQNY
\pk{\sup_{t\in [t_i-\rho,t_i+\rho]} \frac{\chi_{\vk b}^2(t) }{w^2(t)}>u, \sup_{t\in [t_j-\rho,t_j+\rho]} \frac{\chi_{\vk b}^2(t) }{w^2(t)}>u}\leq
\exp\LT(-\frac{(2w(t_1) \sqrt u-Q_{i,j})^2}{2(2+2\eta)}\RT), \quad  1\leq i<j\leq m,
\EQNY
with $Q_{i,j}$'s some positive constants and $\eta\in(0,1)$.
Therefore, as $u\to\IF,$
\BQN\label{doublebound}
\sum_{1\leq i<j\leq m}\pk{\sup_{t\in [t_i-\rho,t_i+\rho]} \frac{\chi_{\vk b}^2(t) }{w^2(t)}>u, \sup_{t\in [t_j-\rho,t_j+\rho]} \frac{\chi_{\vk b}^2(t) }{w^2(t)}>u}=o\left(\MM ( u ) \ \Upsilon_k( w^2(t_1)u)\right).
\EQN
Combining  \eqref{main}--(\ref{doublebound}) with (\ref{decom}) we establish the claim of {\bf C1($\beta$)}.

{\it \underline{\bf C2($\beta$)}}. 
In this case, we have that (\ref{eq:pu}) holds with
\BQNY
\MM_i (\beta_i, u )=\left\{\begin{array}{cc}
 \Piter_{\alpha}^{  a_i(w(t_1)C(t_i))^{-1}\mathcal{L}^\alpha }, & i\in K\\
1,& i\in K^c.
\end{array}
\right.
\EQNY
Consequently,
\BQNY
\sum_{i=1}^m p_i(u)\sim  \left( \prod_{l=k+1}^n (1-b_l^2)^{-1/2}\right)\left(\sum_{i\in K}  \Piter_{\alpha}^{  a_i(w(t_1)C(t_i))^{-1}\mathcal{L}^\alpha }+\sharp K^c\right) \ \Upsilon_k( w^2(t_1)u).
\EQNY
Note that  (\ref{eq:TTr}) and (\ref{doublebound}) still hold. Similarly as the case {\bf C1($\beta$)}, we establish the claim of {\bf C2($\beta$)}.

{\it \underline{\bf C3($\beta$)}}.
In this case, we have that (\ref{eq:pu}) holds with
\BQNY
\MM_i (\beta_i, u )=1, \quad 1\leq i\leq m.
\EQNY
Consequently,
\BQNY
\sum_{i=1}^m p_i(u)\sim  m \left( \prod_{l=k+1}^n (1-b_l^2)^{-1/2}\right) \Upsilon_k( w^2(t_1)u).
\EQNY
Similarly as before, the claim of {\bf C3($\beta$)} follows.

\underline{ (ii).} By {\bf F2} we have  for any sufficiently small $\vn>0$ it holds that
$$
\inf_{t\in\TT_\vn} w(t)>w(c_1),\ \ \text{with} \ \TT_\vn=\TT\setminus\bigcup_{i=1}^m[c_i-\vn,d_i+\vn].
$$
Similarly to \eqref{decom} we have
\BQN\label{eqlower}
&&\sum_{i=1}^m\pk{\sup_{t\in [c_i-\epsilon,d_i+\epsilon]} \frac{\chi_{\vk b}^2(t) }{w^2(t)}>u}+\pk{\sup_{t\in \TT_\vn} \frac{\chi_{\vk b}^2(t) }{w^2(t)}>u}\nonumber\\
&&\ge \pk{\sup_{t\in \TT} \frac{\chi_{\vk b}^2(t) }{w^2(t)}>u}\\
&& \geq\sum_{i=1}^m\pk{\sup_{t\in [c_i,d_i]} \frac{\chi_{\vk b}^2(t) }{w^2(t)}>u}-\sum_{1\leq i<j\leq m}\pk{\sup_{t\in [c_i,d_i]} \frac{\chi_{\vk b}^2(t) }{w^2(t)}>u, \sup_{t\in [c_j,d_j]} \frac{\chi_{\vk b}^2(t) }{w^2(t)}>u}.\nonumber
\EQN
\COM{ 
\BQN
\pk{\sup_{t\in \TT} \frac{\chi_{\vk b}^2(t) }{w^2(t)}>u} \leq  \sum_{i=1}^m\pk{\sup_{t\in [c_i-\epsilon,d_i+\epsilon]} \frac{\chi_{\vk b}^2(t) }{w^2(t)}>u}+\pk{\sup_{t\in \TT_\vn} \frac{\chi_{\vk b}^2(t) }{w^2(t)}>u}.\label{equpper}
\EQN}
Next, we have  from {\bf F2} that for $1\leq i\leq m$
\BQNY
\pk{\sup_{t\in [c_i,d_i]} \frac{\chi_{\vk b}^2(t) }{w^2(t)}>u}&=&\pk{\sup_{t\in [c_i,d_i]}  \chi_{\vk b}^2(t)  >w^2(c_1)u},\\
\pk{\sup_{t\in [c_i-\vn,d_i+\vn]} \frac{\chi_{\vk b}^2(t) }{w^2(t)}>u}&\le&\pk{\sup_{t\in [c_i-\vn,d_i+\vn]}  \chi_{\vk b}^2(t)  >w^2(c_1)u}.
\EQNY
It is noted that   the result in Theorem 2.1 of \cite{LiuJi} also holds when $g(t)=0$. Thus,
it follows from  that result, as $u\to\IF,$
\BQNY
\pk{\sup_{t\in [c_i,d_i]}  \chi_{\vk b}^2(t)  >w^2(c_1)u}&\sim& \prod_{j=k+1}^n\left(1-b_j^2\right)^{-1/2}   \mathcal{H}_{\alpha} \int_{c_i}^{d_i} (C(t))^{1/\alpha}dt  \ (q(w^2(c_1)u) )^{-1}  \Upsilon_k( w^2(c_1)u),\\
\pk{\sup_{t\in [c_i-\vn,d_i+\vn]}  \chi_{\vk b}^2(t)  >w^2(c_1)u}&\sim& \prod_{j=k+1}^n\left(1-b_j^2\right)^{-1/2}   \mathcal{H}_{\alpha} \int_{c_i-\vn}^{d_i+\vn} (C(t))^{1/\alpha}dt  \ (q(w^2(c_1)u) )^{-1}  \Upsilon_k( w^2(c_1)u).
\EQNY
Moreover, since
\BQNY
&& \pk{\sup_{t\in [c_i,d_i]} \frac{\chi_{\vk b}^2(t) }{w^2(t)}>u, \sup_{t\in [c_j,d_j]} \frac{\chi_{\vk b}^2(t) }{w^2(t)}>u}\\
&&\quad = \pk{\sup_{t\in [c_i,d_i]} \chi_{\vk b}^2(t) >w^2(c_1)u, \sup_{t\in [c_j,d_j]} \chi_{\vk b}^2(t) >w^2(c_1)u},
\EQNY
we have from \nelem{LemDouble} that, for all $u$ large,
\BQNY
 \pk{\sup_{t\in [c_i,d_i]} \frac{\chi_{\vk b}^2(t) }{w^2(t)}>u, \sup_{t\in [c_j,d_j]} \frac{\chi_{\vk b}^2(t) }{w^2(t)}>u} \leq
\exp\LT(-\frac{(2w(c_1) \sqrt u-Q_{i,j})^2}{2(2+2\eta)}\RT), \quad  1\leq i<j\leq m,
\EQNY
with $Q_{i,j}$'s some positive constants and $\eta\in(0,1)$.
 This implies that
$$\sum_{1\leq i<j\leq m}\pk{\sup_{t\in [c_i,d_i]} \frac{\chi_{\vk b}^2(t) }{w^2(t)}>u, \sup_{t\in [c_j,d_j]} \frac{\chi_{\vk b}^2(t) }{w^2(t)}>u}=o\left((q(w^2(c_1)u) )^{-1}  \Upsilon_k( w^2(c_1)u)\right), \quad u\rw\IF.$$
Moreover, \nelem{ThmBorell} gives that
$$\pk{\sup_{t\in \TT_\vn} \frac{\chi_{\vk b}^2(t) }{w^2(t)}>u}=o\left((q(w^2(c_1)u) )^{-1}  \Upsilon_k( w^2(c_1)u)\right), \quad u\rw\IF.$$
Consequently, 
by letting $\vn\to0$ we conclude that the claim in (ii) is established.
This  completes the proof. \QED

\proofkorr{corollary} We have from Example \ref{example}, for $\rho_1>0$ and $\rho_2\in \mathbb{R}$,
\BQNY
\sup_{t\in (0,1)}\frac{\chi_{\vk b}^2(t) }{w_{\rho_1, \rho_2}^2(t)}<\IF\ \ a.s.
\EQNY
Furthermore, for the generic locally stationary Gaussian process $X=\overline{B}$ we have
\BQN\label{CK}
C(t)=\frac{1}{2t(1-t)},\ \  \quad K(h)=\sqrt{\abs{h}},\ \ \ h\in(0,1).
\EQN
Next, in order to apply Theorem \ref{asym} we analyze the function   $w_{\rho_1,\rho_2}(t)$.



For simplicity, 
we define
\BQNY 
f(x)=2\rho_1x+2\rho_2\ln x,\ \ \ x(t)=\ln\ln\LT(\frac{e^2}{t(1-t)}\RT).
\EQNY
 Apparently,
$$w_{\rho_1,\rho_2}^2(t)=f(x(t)), \quad t\in (0,1), \quad \{x(t): t\in (0,1)\}=[\ln\ln(4e^2), \IF).$$
Since
$$\frac{\partial f(x)}{\partial x}=2\rho_1+\frac{2\rho_2}{x}, \quad x\in [\ln\ln(4e^2), \IF),$$
the following three different cases will be discussed separately:
\BQNY
a). \ \rho_2> -\rho_1\ln\ln(4e^2); \ \ b). \  \rho_2= -\rho_1\ln\ln(4e^2);\ \  c). \ \rho_2< -\rho_1\ln\ln(4e^2) .
\EQNY

\underline{$a). \ \rho_2> -\rho_1\ln\ln(4e^2)$:} In this case, we have
$$
f'(x)=\frac{\partial f(x)}{\partial x}>0, \quad x\in [\ln\ln(4e^2), \IF),
$$
which means that $f(x)$ attains its minimum over $[\ln\ln(4e^2), \IF)$ at the unique point $x_0=\ln\ln(4e^2)$, and $f'(x_0)>0$. Since further
$$
x'(t)=\frac{\partial x(t)}{\partial t}=\frac{t(1-t)}{\ln\LT(\frac{e^2}{t(1-t)}\RT)} \LT(t^{-1}(1-t)^{-2}-t^{-2}(1-t)^{-1}\RT),
$$
we conclude that the minimizer of $f(x(t))$ over $(0,1)$ is unique and equal to $t_1=1/2$, and $x(t_1)=x_0$,    $x'(t_1)=0$.
Next, we look at the Taylor expansion of  $f(x(t)), t\in(0,1)$ at $t_1$. 
Note that
\BQNY
\frac{\partial f(x(t))}{\partial t}\Big\lvert _{t=t_1}=0, \ \ \ \
\frac{\partial^2 f(x(t))}{\partial t^2}\Big\lvert _{t=t_1}= \frac{\partial f }{\partial x} \Big\lvert _{x=x_0} \frac{\partial^2 x }{\partial t^2}\Big\lvert _{t=t_1}=\frac{ 16 }{ \ln(4e^2)}\LT( \rho_1 +\frac{ \rho_2}{\ln\ln(4e^2)}\RT) =:Q_1>0.
\EQNY
Thus, by Taylor expansion
\BQNY
f(x(t))-f(x(t_1))=\frac{Q_1}{2} (t-t_1)^2\oo,\ \ \ t\to t_1
\EQNY
which yields that
\BQN \label{eq:wa}
w_{\rho_1,\rho_2}(t)=w_{\rho_1,\rho_2}(t_1)+\frac{Q_1}{4 w_{\rho_1,\rho_2} (t_1)}(t-t_1)^2\oo, \quad t\rw t_1.
\EQN
Moreover, by \eqref{CK} and \eqref{eq:wa}
\BQNY
C(t_1)=2,\ \ \alpha=1<2=\beta,\ \  q(u)=u^{-1}.
\EQNY
Since further $\Gamma(1/2+1)=1/2\sqrt{\pi}$ and $\HH_1=1$, by applying \netheo{asym} we conclude that the claim in (a) is established.

\underline{$b). \ \rho_2= -\rho_1\ln\ln(4e^2)$:} In this case, we have that $f(x)$  attains its minimum over $[\ln\ln(4e^2), \IF)$ at the unique point $x_0=\ln\ln(4e^2)$, but with $f'(x_0)=0$. Further,   the minimizer of $f(x(t))$ over $(0,1)$ is unique and equal to $t_1=1/2$, and $x(t_1)=x_0$,    $x'(t_1)=0$.
Thus, we have 
\BQNY
\frac{\partial f(x(t))}{\partial t}\Big\lvert _{t=t_1}=0, \ \ \ \
\frac{\partial^2 f(x(t))}{\partial t^2}\Big\lvert _{t=t_1}=  0.
\EQNY
Next, let us calculate higher-order derivatives of $f(x(t))$. We have
\BQNY 
\frac{\partial^3 f(x(t))}{\partial t^3}&=&\frac{\partial^3 f }{\partial x^3}\LT(\frac{\partial x }{\partial t}\RT)^3 +3\frac{\partial^2 f }{\partial x^2}\frac{\partial x }{\partial t}\frac{\partial^2 x }{\partial t^2} + \frac{\partial f }{\partial x}\frac{\partial^3 x }{\partial t^3},\\
\frac{\partial^4 f(x(t))}{\partial t^4}&=&\frac{\partial^4 f }{\partial x^4}\LT(\frac{\partial x }{\partial t}\RT)^4 +6\frac{\partial^3 f }{\partial x^3}\LT(\frac{\partial x }{\partial t}\RT)^2\frac{\partial^2 x }{\partial t^2} +3\frac{\partial^2 f }{\partial x^2}\LT( \frac{\partial^2 x }{\partial t^2} \RT)^2+4\frac{\partial^2 f }{\partial x^2} \frac{\partial x }{\partial t} \frac{\partial^3 x }{\partial t^3}+ \frac{\partial f }{\partial x}\frac{\partial^4 x }{\partial t^4}.
\EQNY
This implies that
\BQNY
\frac{\partial^3 f(x(t))}{\partial t^3}\Big\lvert _{t=t_1}=0, \ \ \ \
\frac{\partial^4 f(x(t))}{\partial t^4}\Big\lvert _{t=t_1}= 3\frac{\partial^2 f }{\partial x^2}\Big\lvert _{x=x_0}\LT( \frac{\partial^2 x }{\partial t^2}\Big\lvert _{t=t_1} \RT)^2=\frac{384\rho_1}{\ln\ln(4e^2) (  \ln(4e^2))^2} =:Q_2 > 0.
\EQNY
Therefore, by Taylor expansion we conclude that
\BQN\label{eq:wb}
w_{\rho_1,\rho_2}(t)=w_{\rho_1,\rho_2}(t_1)+\frac{Q_2}{48w_{\rho_1,\rho_2} (t_1)}(t-t_1)^4\oo, \quad t\rw t_1.
\EQN
Similarly as in (a), the claim of (b) follows by applying \netheo{asym}.

\underline{$c). \ \rho_2< -\rho_1\ln\ln(4e^2)$:} In this case, we have that $f(x)$ attains its minimum over $[\ln\ln(4e^2), \IF)$ at an inner point $x_0=-{\rho_2}/{\rho_1}$, for which $f'(x_0)=0$. Furthermore, the minimizer of $f(x(t))$ over $(0,1)$ are two distinct points $t_1=\frac{1+\sqrt{1-4e^{2-e^{-\rho_2/\rho_1}}}}{2}$ and $t_2=\frac{1-\sqrt{1-4e^{2-e^{-\rho_2/\rho_1}}}}{2}$, for which $x'(t_i)\neq 0, i=1,2$. Thus, we have, for $i=1,2,$
\BQNY
\frac{\partial f(x(t))}{\partial t}\Big\lvert _{t=t_i}=0, \ \ \ \
\frac{\partial^2 f(x(t))}{\partial t^2}\Big\lvert _{t=t_i}=\frac{\partial^2 f }{\partial x^2} \Big\lvert _{x=x_0}\LT(\frac{\partial x }{\partial t}\Big\lvert _{t=t_i}\RT)^2=\frac{2\rho_1^2}{-\rho_2}Q_{3}>  0,
\EQNY
where, by symmetry of $x(t), t\in(0,1)$,
$$
Q_{3}:= \LT(\frac{\partial x }{\partial t}\Big\lvert _{t=t_1}\RT)^2= \LT(\frac{\partial x }{\partial t}\Big\lvert _{t=t_2}\RT)^2>0.
$$
Consequently, by Taylor expansion we conclude that
\BQN\label{eq:wc}
w_{\rho_1,\rho_2}(t)=w_{\rho_1,\rho_2}(t_i)+\frac{\rho_1^2 Q_3}{-2\rho_2 w_{\rho_1,\rho_2} (t_i)}(t-t_i)^2\oo, \quad t\rw t_i.
\EQN
In addition,
$$
w_{\rho_1,\rho_2} (t_1)=w_{\rho_1,\rho_2} (t_2)=\sqrt{2\rho_2(\ln(-\rho_2) -\ln(\rho_1)-1)}.
$$
Similarly as in (a), the claim of (c) follows by applying \netheo{asym}. This completes the proof. \QED

\proofkorr{corollary1}
First note that
\BQNY
\lim_{h\to 0}\frac{1-\E{\overline{B}_H(t),\overline{B}_H(t+h)}}{ |h|^{2H}}=\frac{1}{2t^{2H}}
\EQNY
holds uniformly in $t\in I$, for any compact interval $I$ in $(0,1]$.
This means that $\overline{B}_H$ is a locally stationary Gaussian process with
\BQNY
C(t)=\frac{1}{2t^{2H}}, \quad K(h)=\abs{h}^{H},\ \ \alpha=2H.
\EQNY
We shall first discuss the finiteness of $\sup_{t\in(0,1]}\frac{\chi_{\vk b}^2(t) }{w_{\rho,\vn}^2(t)}$, for which we only need to verify the conditions in \netheo{Thm01law} for the case where $S=0$.
First, condition {\bf A}(0) is satisfied by  $w_{\rho,\vn}$, and clearly
\BQNY
f(0)=\frac{1}{2^{1/(2H)}} \int_{1/2}^0t^{-1}dt =-\IF.
\EQNY
Further, we have from the calculations in the proof of Corollary 2.7 in \cite{LiuJi} that $\E{\overline{B}_H(t), \overline{B}_H(s)}<1$ for $s\neq t, s,t\in (0,1]$, and conditions {\bf B(0)} and {\bf C(0)} are satisfied by $\overline{B}_H(t), t\in (0,1]$.
Moreover, we have that
\BQNY
J_{c,w_{\rho,\vn}}(0)\leq \rd{\frac{1}{2^{1/(2H)}} \LT( \frac{1}{\LT(\ln\LT(e^2/\epsilon\RT)\RT)^{\rho c}}  \int_{\vn \wedge 1/2}^{1/2}\frac{1}{t }dt+      \int_0^{\vn}\frac{1}{t \LT(\ln\LT(e^2/t\RT)\RT)^{\rho c}}dt \RT)<\IF}
\EQNY
for any $c>1/\rho$.
Consequently, it follows from \netheo{Thm01law}  that
$$
\sup_{t\in(0,1]}\frac{\chi_{\vk b}^2(t) }{w_{\rho,\vn}^2(t)}<\IF\ \ \ \  a.s.
$$
Next, since by definition $w_{\rho,\vn}$ attains its minimum over $(0,1]$ on the interval $[\vn,1]$, we conclude from (ii) of \netheo{asym} that the claim follows. \QED

\bigskip
{\bf Acknowledgement}:  
P. Liu was partially supported by  the Swiss National Science Foundation Grant 200021-166274.

\bibliographystyle{ieeetr}

 \bibliography{GCC1}
\end{document}